\documentclass[a4paper,12pt,twoside]{article} 
\usepackage{german}
\usepackage{fleqn}
\usepackage{eucal}
\usepackage{epsf}
\usepackage[OT1]{fontenc}
\usepackage[utf8]{inputenc}
\usepackage{amsmath}
\usepackage{amsfonts}
\usepackage{amssymb}
\usepackage{bbm}
\usepackage{color}
\usepackage{graphicx}
\usepackage{hyperref}
\usepackage{lscape}
\usepackage{tikz}
\hypersetup{pdftex,pdfstartview=FitH}
\hypersetup{colorlinks=true,urlcolor=blue,citecolor=blue,linkcolor=blue}
\newfont{\qqq}{cminch scaled 1440}
\newfont{\qqqq}{cminch scaled 2074}
\newcommand{\bes}{\begin{eqnarray*}}
\newcommand{\ees}{\end{eqnarray*}}
\newcommand{\bi}{\begin{itemize}\itemsep0em}
\newcommand{\ei}{\end{itemize}}
\newcommand{\gdw}{ \Longleftrightarrow }
\newcommand{\be}[1]{\begin{eqnarray} \label{#1}}
\newcommand{\ee}{\end{eqnarray}}
\newcommand{\darfo}{ \Longrightarrow }
\definecolor{rot}{rgb}{1.000,0.000,0.000}
\definecolor{gruen}{rgb}{0.500,0.86,0.43}

\newcommand{\eyy}{= \\ &=&}
\newcommand{\Z}{{\mathbb{Z}}}
\renewcommand{\P}{{\mathbb{P}}}

\newcommand{\oA}{{\bf A}}
\newcommand{\oI}{{\bf I}}
\newcommand{\oS}{{\bf S}}
\newcommand{\pf}{ \longrightarrow }
\newcommand{\bn}{\begin{enumerate}}
\newcommand{\en}{\end{enumerate}}

\begin{document}
\author{\it Holger Stephan\thanks{e-mail: {\tt stephan@wias-berlin.de}~~~~~~
URL: {\tt http://www.wias-berlin.de/people/stephan}},
Berlin}
\title{\sc \vspace*{-2cm}
  Millions of Perrin pseudoprimes including a few giants}
% in English: lowercase; start with capital letter after colon or dash; short title for page headings
%\nopreprint{2657}	% preprint number
%\nopreyear{2019}	% preprint year
\selectlanguage{english}		% do not change; important for date format
\date{\today}			% fix date, e.g. February 22, 2017
%\date{December 16, 2019}			% fix date, e.g. February 22, 2017
%\subjclass[2010]{11B37, 11B39, 11B50}	% Math. Subject Classif.
%\pacs[2008]{}				% Physics Astronomy Classif., if any
%\keywords{Pseudoprimes, recurrence sequences, fast algorithm, large numbers}
%\thanks{}				% acknowledgements; period is set automatically!
%% amsart only! abstract before maketitle
%\begin{abstract}Abstract\ldots\end{abstract}
\maketitle
%% article and other classes: abstract after maketitle
\begin{abstract}
 The calculation of many and large Perrin pseudoprimes is a challenge.
  This is mainly due to their rarity. Perrin pseudoprimes  are one of the
  rarest known pseudoprimes. In order to calculate many such
  large numbers, one needs not only a fast algorithm but also
  an idea how most of them are structured to minimize the amount of
  numbers one have to test.
  
  We present a quick algorithm for testing Perrin pseudoprimes and
  develop some ideas on how Perrin pseudoprimes might be structured.
  This leads to some conjectures that still need to be proved.
  
  We think that we have found well over 90\% of all 20-digit Perrin
  pseudoprimes. Overall, we have been able to calculate over 9 million
  Perrin pseudoprimes with our method, including some very large ones.
  The largest number found has 3101 digits.
  This seems to be a breakthrough, compared to the previously known just
  over 100,000 Perrin pseudoprimes, of which the largest have 20 digits.

  In addition, we propose two new sequences that do not provide any
  pseudoprimes up to $10^9$ at all.

  \end{abstract}
%%%%%%%%%%%%%%%%%%%%%%%%%%%%%%%%%%%%%%%%%%%%%%%%%%%%%%%%%%%%%%%%%%%%%%
%                                                                    %
%                         Start the article here                     %
%                                                                    %
%%%%%%%%%%%%%%%%%%%%%%%%%%%%%%%%%%%%%%%%%%%%%%%%%%%%%%%%%%%%%%%%%%%%%%

\tableofcontents

%----------------------------------------------------------------------
%\newpage
%~\newpage
%\vspace*{-1cm}

\section{Introduction}

To motivate that it makes sense to deal with primes,
it is best to quote Gauss \cite{gauss}:

~

\begin{center}
  \begin{minipage}{15cm}
\begin{center}
  {\it
    The problem of distinguishing prime numbers from composite numbers, and of
    resolving the latter into their prime factors is known to be one of the
    most important and useful in arithmetic. It has engaged the industry
    and wisdom of ancient and modern geometers to such an extent that it
    would be superfluous to discuss the problem at length. Nevertheless
    we must confess that all methods that have been proposed thus
    far are either restricted to very special cases or are so laborious
    and difficult that even for numbers that do not exceed the limits
    of tables constructed by estimable men, they try the patience of
    even the practiced calculator. And these methods do not apply at
    all to larger numbers.
}
\end{center}
  \end{minipage}
\end{center}

~
  
Prime numbers are a very serious issue.
We prefer dealing with pseudoprimes.
Pseudoprimes are numbers that behave similar to primes.
%Unfortunately, we are not aware of a quote from Gauss to pseudoprimes.
%But in Gauss's quote we can replace the word prime with pseudoprime.

Sometimes it is a big challenge to compute all or at least many
or some very large pseudoprimes of a given type.

In this paper, we introduce a quick algorithm for the calculation
of Perrin pseudoprimes. This is nothing special,
there are already many fast algorithms.
Similar to primes, also for pseudoprimes it is difficult to guess their
structure. Therefore, in order to calculate all of them there is nothing left
but to test every single number. This strongly limits the size of the numbers.
It turns out, however, that the structure can be guessed for
most of the pseudoprimes.
This very much limits the range of potential numbers to be tested
and makes it
possible to calculate millions of them and even very large ones.

We do the following {\bf Notations:}
\bi
\item The set of all primes is denoted by $\P$.
\item $a|b$ means $a$ divides $b$ or $b$ is divisible by $a$.
\item We state some classical facts from number theory as theorems,
  omitting the proofs.
\ei

\subsection{The Perrin sequence}

Let us define a sequence (called Perrin sequence) $P_n$ recursively:
\bes
P_0 &=& 3\\
P_1 &=& 0\\
P_2 &=& 2 \\
P_n &=& P_{n-2} + P_{n-3} ,~n \geq 3
\ees
and calculate the first entries:
\bes
(P_n)_{n=0}^\infty  = 3, 0, 2, 3, 2, 5, 5, 7, 10, 12, 17, 22, 29, 39, 51, 68, 90, 119, ...
\ees
We observe: If $n$ is prime, we have $n| P_n$ and that goes on for a long time.

Anyone seeing this sequence for the first time is certainly quite surprised,
since it is believed that there is no simple algorithm for calculating
the primes.

The recursion law of this sequence was found in 1899 by Edouard Lucas.
This sequence with the initial values given above, was first used by
Raoul Perrin \cite{wiki1,oeis}.

Probably many mathematicians and amateur mathematicians have tried to
answer the question of whether this sequence really only produces primes.
Considering that already the number $P_{811}$ has 100 digits,
one can imagine how difficult that has been.

The answer was not found until 1982, when Jeffrey Shallit
(according to \cite{oeis}) calculated  the first
two non-prime numbers -- so-called Perrin pseudoprimes (PPP) --
with a computer.
Here they are: $ 271441 = 521 \cdot $ 521 and 
$ 904631 = 7 \cdot 13 \cdot 9941 $. $ P_{271441} $ has 33150 digits.
Today it is known that there are infinitely many
Perrin pseudoprimes \cite{granv1}. Nevertheless,
they are very rare, which makes their finding still difficult.

In this paper, we develop an effective algorithm for calculating
Perrin pseudoprimes and present some numerical results that
constitute, to our knoledge, right now
the world's largest collection of Perrin
pseudoprimes including the largest PPP.

%https://en.wikipedia.org/wiki/Perrin_number#Perrin_pseudoprimes

\section{Pseudoprimes}

\subsection{Iff-- and if--Theorems}

There are two kinds of theorems dealing with primes that can be used
to test a given number on whether it is a prime.

1) Theorems like: $p \in \P$ if and only if property $A(p)$ holds.

2) Theorems like: $p \in \P$, then property $A(p)$ holds.

Theorems of the first kind are, for example
\bi
\item {\bf Theorem:} $p \in \P ~\gdw~ \forall k\in \P, k \leq \sqrt{p}:~k\not| p$
\item {\bf Theorem (Wilson):} $
p \in \P ~\gdw~ p| 1\cdot 2  \cdot 3 \cdots (p-1) + 1$
\item {\bf Theorem:}
  \be{e201}
  p \in \P ~\gdw~ p|{p \choose k} ~~~\forall k=1,...,p-1
  \ee
\ei
These theorems allow for deterministic tests. If for a given number $p$
the property $A(p)$ holds, then $p$ is prime.

Unfortunately, algorithms based on deterministic testing have
high complexity, so far.

~

Theorems of the second kind state:  If for a given number $p$
the property $A(p)$ holds, then $p$ can be prime or not. This is
useful, if $p$ is prime with very high ``probability''. Testing
$A(p)$ one can be ``very sure'' that $p$ is prime.
Typically such kind of probabilistic tests are much faster (have
a lower complexity) than deterministic ones. 
Thus, it is useful to create tests with a very small
equivalence gap, the gap between if and iff.

Numbers $n$ that lie in this gap, i.e. $A(n)$ holds, but
$n$ is composite, are called pseudoprimes with respect to  property $A$.

One example, following immediately from (\ref{e201}) is:

{\bf Theorem:}
 $p \in \P  ~~\darfo~~ 
p\left|{\displaystyle \sum_{k=1}^{p-1} a_k {p \choose k}} \right.$
for some given integers $a_k$.

~

It is clear that looking at a linear combination of binomial coefficients
instead of all coefficients in detail, we loose information. This is
just the equivalence gap.
Looking at a given linear combination  of binomial coefficients is
faster than looking at every one in detail. The idea is to
choose such coefficients $a_k$ so that the equivalence gap is small.

~

Here, we define some kind of probability (better frequency)
for a pseudoprime test.
Let $\pi(n)$ be the number of primes less than $n$ and $P(n)$ the 
number of pseudoprimes less than $n$ for a given pseudoprime test. By
$W(n) = P(n) / \pi(n)$ we define the frequency of numbers incorrectly tested
and call it error rate.
Thus, the lower the error rate $W(n)$, the better the test.

Of course, it would be best if a test provided
only a finite number of pseudoprimes.
These would be calculated and stored in a database which allowed for 
a deterministic test, practically. 
Such a test is not yet known. In contrast, until now,
for many pseudoprime number type, it has been proved sooner or later
that there are infinitely many ones.

%\newpage

\subsection{Fermat and Carmichael pseudoprimes}

The simplest pseudoprimes are Fermat pseudoprimes. They are consequences
of Fermat's little

{\bf Theorem:} Given an integer $z \geq 2$.
If $p \in\P$ then $\displaystyle{p| z^p - z} $.

Conversely, if a number $n$ for some $z$ satisfies 
$\displaystyle{n| z^n - z} $ but $n \not\in\P$, $n$ is called
Fermat$_z$ pseudoprime.

Best known is the special case $z=2$:

{\bf Theorem:} If $p \in\P$ then $\displaystyle{p| 2^p - 2} $.

A number $n \not\in\P$ with $\displaystyle{n| 2^n - 2} $
is called Fermat$_2$ pseudoprime.

\subsubsection{Fermat$_2$ pseudoprimes}

Fermat's little Theorem for $z=2$ is an easy consequence of
Theorem 1.

Indeed, multiplying out $(a+b)^n$ with integers $a,b$ we get
\bes
(a+b)^n &=& a^n +
 {n \choose 1} a^{n-1} b +  {n \choose 2} a^{n-2} b^2 + 
{n \choose 3} a^{n-3} b^3 +  \cdots  + b^n
\ees
Therefore, defining
\bes
f_n = (a+b)^n - a^n -b^n =  {n \choose 1} a^{n-1} b + \ldots +
{n \choose n-1} a b^{n-1}
\ees
we obtain the

{\bf Theorem:} If $p \in\P$ then $\displaystyle{p|f_p} $.

The special case ($a=b=1$) yields Fermat's little Theorem to the base $z=2$.

~

Let's calculate the first ones:

~

\begin{tabular}{|r|r|r|r|r|}
\hline
$n$  & $2^n - 2$ & $n|2^n - 2$ ? & $n$ is prime? \\
\hline
2 & 2 & {\textcolor{gruen}{yes!}} & {\textcolor{gruen}{yes!}}  \\
3 & 6 & {\textcolor{gruen}{yes!}}  & {\textcolor{gruen}{yes!}}  \\
4 & 14 & no! & no! \\
5 & 30 & {\textcolor{gruen}{yes!}}  & {\textcolor{gruen}{yes!}}  \\
6 & 62 & no! & no! \\
7 & 126 & {\textcolor{gruen}{yes!}}  & {\textcolor{gruen}{yes!}}  \\ %\pause
\hline
341 & 4479... (103 digits) & {\textcolor{gruen}{yes!}}   &  
{\textcolor{rot}{no!}}~~ $341=11 \cdot 31$ ~~~~\\
561 & 7547... (169 digits) & {\textcolor{gruen}{yes!}}   & 
{\textcolor{rot}{no!}}~~ $561=3 \cdot 11 \cdot 17$\\
645 & 1459... (195 digits) & {\textcolor{gruen}{yes!}}   & 
{\textcolor{rot}{no!}}~~ $645= 3\cdot 5  \cdot 43$ ~~\\
\hline
\end{tabular}

~

Up to 100000 we have 78 pseudoprimes and 9592 primes. 
Thus, we have $W(10^5) = 0.00813178$.

\subsubsection{Carmichael numbers}

Instead of $z=2$ we can consider Fermat$_z$ pseudoprimes with other bases.
Maybe other bases provides fewer pseudoprimes? It turns out that
$z=2$ is one of the best bases. Moreover, there are 
non-primes $n$ with $n|z^n-z$ for any base $z$, the so-called
Carmichael numbers. 561 is the smallest one. The next ones are

~

\begin{minipage}{5cm}
\begin{tabular}{|c|c|}
\hline
Carmichael number & factors \\ \hline
561 &    3 $\cdot$ 11 $\cdot$ 17 \\
1105  &		5 $\cdot$ 13 $\cdot$ 17\\
1729  &		7 $\cdot$ 13 $\cdot$ 19\\
2465  &		5 $\cdot$ 17 $\cdot$ 29\\
2821  &		7 $\cdot$ 13 $\cdot$ 31\\
6601  &		7 $\cdot$ 23 $\cdot$ 41\\
8911  &		7 $\cdot$ 19 $\cdot$ 67\\
10585  &	5 $\cdot$ 29 $\cdot$ 73\\
\hline
\end{tabular}
\end{minipage}
~~~~~~~~~~~~~~~~~~~~~~~
\begin{minipage}{5cm}
\begin{tabular}{|c|c|}
\hline
Carmichael number & factors \\ \hline
15841  &		7 $\cdot$ 31 $\cdot$ 73\\
29341  &		13 $\cdot$ 37 $\cdot$ 61\\
41041  &		7 $\cdot$ 11 $\cdot$ 13 $\cdot$ 41\\
46657  &		13 $\cdot$ 37 $\cdot$ 97\\
52633  &		7 $\cdot$ 73 $\cdot$ 103\\
62745  &		3 $\cdot$ 5 $\cdot$ 47 $\cdot$ 89\\
63973  &		7 $\cdot$ 13 $\cdot$ 19 $\cdot$ 37\\
75361  &		11 $\cdot$ 13 $\cdot$ 17 $\cdot$ 31\\
\hline
\end{tabular}
\end{minipage}

~

There are 16 Carmichael numbers up to 100000.
Moreover, we have the following

{\bf Theorem:} There are infinitely many Carmichael numbers \cite{granv1}.

\subsection{General pseudoprimes} 

\subsubsection{Sums of powers. Multinomial coefficients}

Similar to binomial coefficients, there is a theorem for
multinomial coefficients:

{\bf Theorem:}  $p \in \P  ~~\gdw~~ p|{ p! \over i! ~ j! ~ k!},~\forall
i,j,k \mbox{ with } 0 < i+j+k = p$

~

From this, we conclude the following

{\bf Theorem:}  $p \in \P  ~~\darfo~~ 
p\left|{\sum_{0 < i+j+k=p}~~a_{ijk}~~ { p! \over i! ~ j! ~ k!} }\right.$.
for some integer coefficients $a_{ijk}$. 

From this, multiplying out $(a+b+c)^n$ with integers $a,b,c$
we conclude the

{\bf Theorem:} Given a sequence
\bes
f_n = (a+b+c)^n - a^n -b^n - c^n ~~=~~
\sum_{0 < i+j+k=n} 
{ n! \over i! ~ j! ~ k!} ~ a^i b^j c^k
\ees
Then, $p \in\P$ implies $\displaystyle{p|f_p} $.

~

Similarly we get the

{\bf Theorem:} Given integers $a_1,...,a_k$. Build the sequence
\be{e111}
f_n = (a_1+a_2+ ... + a_k)^n - (a_1^n+a_2^n+ ... + a_k^n)
\ee
Then, $p \in\P$ implies $\displaystyle{p|f_p} $.

~

The example $a_i=1$ yields $f_n = k^n - k$,
Fermat's little theorem in the general case.

~

Perrin's sequence is given in a recurrent way. Here, we recall the
important connection between polynomials and recurrence sequences.

\subsubsection{Polynomials and recurrence sequences}

A linear recurrence sequence (or linear difference equation) of order $k$
is a sequence $(h_n)_{n=0}^\infty$ defined in the following way:

Given $k$ numbers $c_1,...,c_k$ set
\be{e101}
h_n = c_1 h_{n-1} + c_2 h_{n-2} + ... +  c_k h_{n-k} ~.
\ee
Together with $k$ initial conditions $h_0$, $h_1$, ..., $h_{k-1}$
such a sequence is uniquely determined.

~

Obviously, if $c_1,...,c_k$ and $h_0$, $h_1$, ..., $h_{k-1}$,
are integers, then $h_n$ is an integer for all $n$.

~

There is a remarkable connection between such sequences and
polynomials of degree $k$.
If we put $h_n=x^{n}$ and multiply by $x^{k-n}$, we get
an algebraic equation for the roots of
a polynomial formed from the coefficients of the sequence
\be{e103}
Q(x) = - x^k + c_1 x^{k-1} +  c_2 x^{k-2} +  ... + c_{k-1} x + c_k~.
\ee
This polynomial has $k$ -- in general complex -- roots
$x_1,...,x_k$. For simplicity, we assume that the roots are different.

Set
\be{e102}
g_n = b_1 x_1^n + b_2 x_2^n + ... + b_k x_k^n~,
\ee
with some coefficients $b_1,...,b_k$. Solve the system of $k$ linear equations
$h_i=g_i$, $i=0,...,k-1$ with respect to the unknown $b_j$. This 
is always uniquely solvable, because the 
corresponding matrix is the Vandermonde matrix $(x_i^j)$. Its 
determinant does not vanish if the roots $x_i$ are different, as required.

{\bf Theorem:} For any $n \geq 0$ we have $g_n=h_n$.

This is easily proved, since we have $Q(x_i) = 0$ for $i=1,...,k$.

~

The opposite is also true:

{\bf Theorem:} Given $k$ different complex numbers $x_1,...,x_k$ and $k$
real numbers $b_1,...,b_k$. Calculate the first entries $h_0, ..., h_{k-1}$
of some sequence $(h_n)$
by the right-hand side of (\ref{e102}) and compile a
polynomial (\ref{e103}) from it's roots $x_1,...,x_k$
\bes
Q(x) = -(x-x_1) \cdots  (x-x_k) = -x^k + (-1)^{k+1}(x_1 + ... + x_k) x^{k-1} +
...
\ees
Then, the sequence (\ref{e101}), given in a recurrent way is 
exactly the sequence (\ref{e102}), given explicitely.

~

Thus, we have a one-to-one correspondence between the
linear recurrence sequence (\ref{e101}) and the sum of powers (\ref{e102}).

This can be applied to Perrin's sequence.

\subsubsection{Perrin's sequence, given explicitely}

Starting with the sequence
\bes
P_0 &=& 3\\
P_1 &=& 0\\
P_2 &=& 2 \\
P_n &=& P_{n-2} + P_{n-3} ,~n \geq 3
\ees
at first, we  compile the polynomial from the coefficients
\bes
Q(x) = -x^3 + x + 1
\ees
Its roots are
\bes
a &=& 1.32472...\\
b &=& -0.662359... + 0.56228... i \\
c &=& -0.662359... - 0.56228... i
\ees
Set $h_n = a^n + b^n + c^n$ (since $a+b+c = 0$). The
first entries are
\bes
h_0 &=& a^0 + b^0 + c^0 = 3\\
h_1 &=& a^1 + b^1 + c^1 = 0\\
h_2 &=& a^2 + b^2 + c^2 = (a+b+c)^2 -2(ab+bc+ca) = 0 -2(-1) = 2 
\ees
Thus, the sequences $P_n$ and $h_n$ coincide.

~

The theorem
\bes
p \in \P ~~\darfo~~ p|P_p = a^p + b^p + c^p 
\ees
does not follow from this, immediately, since $a,b,c$ are not integers.

We have to answer two questions:
\bi
\item When is $f_n = (a+b+c)^n - a^n -b^n - c^n$ an integer sequence?
\item When does $(p \in\P~\darfo~\displaystyle{p|f_p} )$ hold?
\ei

\subsubsection{When is $f_n = (a+b+c)^n - a^n -b^n - c^n $ an integer?}

For any $n$, the expression $  (a+b+c)^n - a^n -b^n - c^n$
is a symmetric polynomial in $a$, $b$ and $c$.

{\bf Theorem:}
Any symmetric polynomial can be expressed in terms of elementary symmetric 
polynomials.

Here, these are
\bes
A_1=a+b+c,~A_2=ab+bc+ca,~A_3=abc
\ees
which are the coefficients of a polynomial with roots $a,b,c$.

Calculating, for example, the first entries, we get
\bes
(a+b+c)^0 - a^0 - b^0 - c^0 &=& -2\\
(a+b+c)^1 - a^1 - b^1 - c^1 &=& 0\\
(a+b+c)^2 - a^2 - b^2 - c^2 &=& 2 A_2\\
(a+b+c)^3 - a^3 - b^3 - c^3 &=& 3 A_1 A_2 - 3 A_3 \\
(a+b+c)^4 - a^4 - b^4 - c^4 &=& 4 A_1^2 A_2 - 4 A_1 A_3 -  2 A_2^2
\ees

Hence, $f_n$ is integer if $a,b,c$ are roots of a 
polynomial with integer coefficients.

\subsubsection{When does $(p \in\P~\darfo~\displaystyle{p|f_p} )$ hold?}

We have
\bes
f_n=(a+b+c)^n -(a^n + b^n + c^n) = \sum_{0 < i+j+k=n} 
{ n! \over i! ~ j! ~ k!} ~ a^i b^j c^k
\ees
and $p \in \P  ~~\darfo~~ p|{ p! \over i! ~ j! ~ k!},~\forall
i,j,k \mbox{ with } 0 < i+j+k = n$.

$\displaystyle { n! \over i! ~ j! ~ k!}$ does not change by a permutation
of $i,j,k$. It can be lifted out.
\bes
 \sum_{0 < i+j+k=n} 
{ n! \over i! ~ j! ~ k!} ~ a^i b^j c^k
=
 \sum_{0 < i \leq j \leq k} 
{ n! \over i! ~ j! ~ k!} \sum_{\pi(i,j,k)} a^i b^j c^k
\ees
$\sum_{\pi(i,j,k)} a^i b^j c^k$ is again a symmetric polynomial and so
it is an integer if $a,b,c$ are roots of an  polynomial with integer
coefficients.

Hence, if $a,b,c$ are roots of a polynomial with integer coefficients,
and $f_n=(a+b+c)^n -(a^n + b^n + c^n) $, then
$p \in\P~\darfo~\displaystyle{p|f_p} $.

\subsubsection{The recurrent calculation of the sequence}

From the polynomial $Q(x)$ it is easy to compile the 
recurrent relation
\bes
g_{n} =  a_1~ f_{n-1} + a_2~ f_{n-2} + a_3~ f_{k-3} + \ldots + a_k~ f_{n-k}
\ees
corresponding to the explicit expression
\bes
g_{n} = x_1^n + \ldots + x_k^n~.
\ees
From this  explicit expression we have to calculate the initial values
$g_0,...,g_{k-1}$. Then, we have
\bes
f_n = g_n - a_1^n~.
\ees
Actually, this is practicable if $a_1 = 0$ (like in the Perrin case)
or $a_1 = \pm 1$.  In other cases, $a_1^n$ increases rapidly and it is better
to look on
\bes
f_n  = (x_1^n + \ldots + x_k^n) - (x_1 + \ldots + x_k)^n
\ees
as on a sum of $k+1$ powers. This corresponds
to a sequence of order $k+1$, having a corresponding polynomial with 
the $k+1$ roots $x_1,...,x_k, a_1 = x_1 + \ldots + x_k$. This
polynomial is
\bes
G(x) &=& -(x-x_1)\cdots(x-x_1)(x-x_1-\ldots - x_k) =
Q(x) (x-a_1)
\eyy
-x^{k+1} + 2 a_1 x^k + \sum_{i=1}^{k-1} (a_{i+1} - a_1 a_i) x^{k-i}
  - a_1 a_k
\ees

\subsubsection{The main theorem}

Connecting the last facts together, we finally obtain the

{\bf Main Theorem:} Given a polynomial of degree $k$
\bes
Q(x) = -x^k + a_1 ~x^{k-1} + a_2~ x^{k-2} + a_3~ x^{k-3} +
\ldots + a_{k-1}~ x + a_k
\ees
with integer coefficients $a_i \in \Z$ and (maybe complex)
roots $x_1,...,x_k$. Then, the sequence
\bes
f_n = (x_1^n + \ldots + x_k^n) - (x_1 + \ldots + x_k)^n
\ees
is an integer sequence and it holds ~~~$p \in \P ~~\darfo~~ p|f_p$.

The sequence $f_n$ can be calculated in a recurrent way from an
order $k$-recurrent relation 
\bes
g_{n} =  a_1~ f_{n-1} + a_2~ f_{n-2} + a_3~ f_{k-3} + \ldots + a_k~ f_{n-k}
\ees
by $f_n = g_n -  a_1^n$ or directly from an order $(k+1)$-recurrent relation 
\bes
f_n =  2 a_1 f_{n-1} + \sum_{i=1}^{k-1} (a_{i+1} - a_1 a_i) f_{n-i-1} 
- a_1 a_k f_{n-k-1} 
\ees

~

We can conclude that any polynomial with integer coefficients is cantidate
to generate pseudoprimes.

~

%\rb{Analogen Satz für Gaußsche integer?}

\newpage

\section{Numerical algorithms}

To calculate pseudoprimes, at first we have 
to calculate $f_n$ by a recurrent or explicit expression and then
we test whether $n|f_n$.

The recurrence relation seems to be very fast, with some additions
for every number. Unfortunately, the entries $f_n$ grow very fast.
For the Perrin sequence we have $P_n \sim 1.32472...^n$ (the largest root).
Thus, $P_{271441}$ has 33150 decimal digits,
$P_{99607901521441}$ -- the 17-th Perrin pseudoprime has
12,164,524,642,561 decimal digits requiring $\sim$ 5 TByte to store it.

The same problem arises with the explicit expression. We have to calculate
$x_j^n$ considering a huge number of digits to get an integer in the end.
But this is necessary to check the remainder of $f_n$ when divided by $n$.

~

The only useful method is to carry out all operations modulo $n$.
This will save us from the usage of the huge numbers $f_n$. We can
still use the recurrence relation but for every new
number we have to start at the very beginning of the sequence, since
calculating $f_n$ mod $n$, we cannot use the result to
calculate $f_{n+1}$ mod $n+1$.

Even doing so, this is still a problem if we want (and we want!) to deal
with large numbers $n$ having, say, 100 digits.
Note, this is the number of digits of the index, not of the sequence member!

Thus, if $n=10^{100}$ we need a fast algorithm for $10^{100}$ additions
of numbers like $10^{100}$ (all done modulo $n$).

Clearly, this has to be an algorithm with logarithmic complexity.
This can be done in pursuing following steps:
\bi
\item We can calculate $k$ entries of the sequence at once, using
  matrix powers. 
\item The $n$-th power of a matrix can be performed in $\log_bn$
  operations using the decomposition of $n$ with respect to
  a fixed basis and Horner's method.
\item In some special cases -- and the Perrin sequence is such a case --
  the calculation can be further simplified. 
\ei

\subsection{Matrix powers instead of additions}

Given a recurrence sequence of order $k$
\be{e131}
f_n = c_{k-1} f_{n-1} + c_{k-2} f_{n-2} + ... +  c_0 f_{n-k} 
\ee
with initial values
\be{e132}
F_0 := (f_0,...,f_{k-1})~.
\ee
The $k$-th entry
\bes
f_k &=& c_{k-1} f_{k-1} + c_{k-2} f_{k-2} + ... +  c_0 f_{0} 
\ees
is a linear combination of the initial values and so are all entries,
for example the $k+1$-th entry
\bes
f_{k+1} &=& c_{k-1} f_{k} + c_{k-2} f_{k-1} + ... +  c_0 f_{1}
\eyy
c_{k-1} (c_{k-1} f_{k-1} + c_{k-2} f_{k-2} + ... +  c_0 f_{0})
+ c_{k-2} f_{k-1} + c_{k-3} f_{k-2} + ... +  c_0 f_{1}
\eyy
( c_{k-1}^2 + c_{k-2} ) f_{k-1} + (c_{k-1} c_{k-2} + c_{k-3}) f_{k-2} + ... +
(c_{k-1} c_1 +  c_0 ) f_1 +
c_{k-1} c_0 f_{0}
\ees
Writing all the entries $F_1 := (f_k,...,f_{2k-1})$ as linear combinations
of $F_0 = (f_0,...,f_{k-1})$, we can compile
a matrix $\oA$ and write $F_1 = \oA F_0$, i.e.,
\bes
\left(
  \begin{array}{c}
    f_{k} \\ f_{k+1} \\ \vdots \\ f_{2k-1}
\end{array}
\right)
=
\left(
\begin{array}{cccc}
 c_0 & c_1 & \cdots & c_{k-1} \\
  c_{k-1} c_0 & c_{k-1} c_1 +  c_0 & \cdots &  c_{k-1}^2 + c_{k-2} \\
  \vdots  & \vdots  & \vdots  & \vdots  
\end{array}
\right)
\left(
  \begin{array}{c}
    f_{0} \\ f_{1} \\ \vdots \\ f_{k-1}
\end{array}
\right)
\ees
This is an equivalent description of (\ref{e131}),  (\ref{e132}).

In the special case $k=3$ we have
\bes
\left(
  \begin{array}{c}
    f_{3} \\ f_{4} \\ f_5
\end{array}
\right)
=
\left(
\begin{array}{ccc}
 c_0 & c_1 & c_2 \\
 c_0 c_2 & c_0+c_1 c_2 & c_2^2+c_1 \\
 c_0 c_2^2+c_0 c_1 & c_1^2+c_2^2
   c_1+c_0 c_2 & c_2^3+2 c_1 c_2+c_0 \\
\end{array}
\right)
\left(
  \begin{array}{c}
    f_{0} \\ f_{1} \\ f_2
\end{array}
\right)
\ees
It follows $F_m = \oA^m F_0$ for
$F_m = (f_{mk},f_{mk+1},..., f_{(m+1)k-1})$.
Thus, if we want to know $f_n$, we have to divide $n$ by $k$ with
remainder, i.e., to write $n=mk+i$ with $i=0,...,k-1$
and calculate $\oA^m$. Instead of additions we have to calculate
the power of a matrix. This can be done very effectively.

\subsection{Horner's method instead of matrix powers}

We have to calculate $\oA^m$ for a given matrix $\oA$.
Let $m=a_0 b^j + ... + a_{j-1} b + a_j$
be the decomposition of $m$ to base $b$ with $a_0 > 0$
and $b > a_i \geq 0$.
Then, calculating the polynomial $a_0 b^j + ... + a_{j-1} b + a_j$
with Horner's method, iteratively
\bes
a_0 b^j + ... + a_{j-1} b + a_j =
((( (a_0 b) b + a_1 ) b + a_2 )  b + ... + a_j ) 
\ees
we conclude
\bes
\oA^m = \oA^{a_0 b^j + ... + a_{j-1} b + a_j} =
(((((\oI^b \oA^{a_0})^b \oA^{a_1})^b \oA^{a_2})^b \cdots  )\oA^{a_j})
\ees
The vector 
\bes
(\oA_0,\oA_1, ..., \oA_{b-1} ) = (\oA^0,\oA^1,\oA^2,...,\oA^{b-1} )
\ees
can be calculated and stored in advance. The calculation runs
especially effectively if $b$ itself is a power of 2.
For practicle purposes we used $b=2,4,8$.

\subsection{A fast algorithm for the Perrin sequence}

The following algorithm was written in 1982 by Frank Bauern\"oppel and
Uwe Kaufmann \cite{alg} in Berlin.

{\bf 1st step:} 
Given $n$. Set $n=3m+i$, $i\in\{0,1,2\}$.
Since we have
\bes
P_3 &=& P_1 + P_0 \\
P_4 &=& P_2 + P_1 \\
P_5 &=& P_3 + P_2 =  P_1 + P_0 + P_2
\ees
we can introduce a matrix
\bes
\oS = 
\left(
\begin{array}{ccc}
1 & 1 & 0 \\
0 & 1 & 1 \\
1 & 1 & 1 \\
\end{array}
\right)
\ees
and have
\bes
\left(
\begin{array}{c}
P_{3m} \\ P_{3m+1} \\ P_{3m+2}
\end{array}
\right)
&=&
\oS^m
\left(
\begin{array}{c}
3 \\ 0 \\ 2 
\end{array}
\right)
=
\left(
\begin{array}{ccc}
1 & 1 & 0 \\
0 & 1 & 1 \\
1 & 1 & 1 \\
\end{array}
\right)^m
\cdot
\left(
\begin{array}{c}
3 \\ 0 \\ 2 
\end{array}
\right)
\ees

{\bf 2nd step:}
The power of $\oS$ can be further simplified by using the square
$\oS^2$. Depending on whether $m$ is even or odd, one have
\bes
\oS^m = \big( \oS^{m \over 2} \big)^2,~~~2|m ~~~~\mbox{ or }~~~~
\oS^m = \big( \oS^{m-1 \over 2} \big)^2 \cdot \oS,~~~2\not|m
\ees
The total power $\oS^m $ can now be calculated iteratively
by using the binary representation of $m$. Let 
$m = (\alpha_0, \alpha_1,  \alpha_2, ...,  \alpha_k, ...)$, $\alpha_0=1$
be the dual number representation of $m$.
We calculate iteratively matrices $\oS_k$
in the following way:
\bes
\oS_{0} &=& \oI\\
\oS_{k+1} &=& 
\left\{
\begin{array}{lrr}
\oS_k^2 & \mbox{if} &\alpha_k = 0\\
\oS_k^2 \cdot \oS & \mbox{if} & \alpha_k = 1\\
\end{array}
\right.
\ees
Then, $\oS^m = \oS_{k_0}$ for some $k_0 < m$.

For example, we have
\bes
\oS^{22} &=& \oS^{10110_{2}} = 
(((( \oI^2 \cdot \oS )^2 \cdot \oS )^2 \cdot \oS )^2 )^2 = \oS^{22}
\ees
For every 0 
(the even digits) one has to square (operation $Q$), for every 1
(the odd digits) one has to square and then to multiply
(operation $QM$).

{\bf 3rd step:} 
Observe that
\bes
\oS^m = 
\left(
\begin{array}{ccc}
1 & 1 & 0 \\
0 & 1 & 1 \\
1 & 1 & 1 \\
\end{array}
\right)^m
=
\left(
\begin{array}{ccc}
a & c & b \\
b & a+b & c \\
c & b+c & a+b \\
\end{array}
\right)
\ees
Thus, one only has to remember the first column $(a,b,c) $ and to know
how this column changes when multiplying $M$ and squaring $Q$.

Operation multiplying $M$:
\bes
\left(
\begin{array}{ccc}
a & c & b \\
b & a+b & c \\
c & b+c & a+b \\
\end{array}
\right)
\left(
\begin{array}{ccc}
1 & 1 & 0 \\
0 & 1 & 1 \\
1 & 1 & 1 \\
\end{array}
\right)
=
\left(
\begin{array}{ccc}
a + b & a + b + c & b + c \\
b + c & a + 2 b + c & a + b + c\\
a + b + c &   a + 2 b + 2 c & a + 2 b + c \\
\end{array}
\right)
\ees
Thus, $M: (a,b,c) \pf  (a+b,b+c,a+b+c)$.

Operation squaring $Q$:
\bes
\left(
\begin{array}{ccc}
a & c & b \\
b & a+b & c \\
c & b+c & a+b \\
\end{array}
\right)^2
\!\!\!
=\!\!
\left(
\begin{array}{ccc}
a^2 \!+\! 2 b c & b^2 \!+\! 2 a c \!+\! 2 b c & 2 a b \!+\! b^2 \!+\! c^2 \\
2 a b \!+\! b^2 \!+\! c^2 & a^2 \!+\! 2 a b \!+\! b^2 \!+\! 2 b c \!+\! c^2 & 
b^2 \!+\! 2 a c \!+\! 2 b c \\
b^2 \!+\! 2 a c \!+\! 2 b c & 2 a b \!+\! 2 b^2 \!+\! 2 a c \!+\! 2 b c \!+\! c^2 & 
a^2 \!+\! 2 a b \!+\! b^2 \!+\! 2 b c \!+\! c^2 \\
\end{array}
\right)
\ees
Thus, $Q: (a,b,c) \pf  (a^2 + 2 b c , b^2 + c^2 + 2ab, b^2 + 2ac + 2bc )$.

\label{bauer}
Furthermore, some numbers $n$ can be excluded from the beginning, 
because we have %\rb{(Warum?)}
\bes
n \equiv 0 \mbox{ mod } 4~\darfo~ f_n \not\equiv 0 \mbox{ mod } 4
~\darfo~ f_n \not\equiv 0 \mbox{ mod } n~.
\ees
The same happens for $n=9,14,...$. Moreover, we have
\bes
n \equiv 0 \mbox{ mod } 3~,~~
n \not\equiv 0,1,3,9 \mbox{ mod } 13
~\darfo~ f_n \not\equiv 0 \mbox{ mod } 3
~\darfo~ f_n \not\equiv 0 \mbox{ mod } n
\ees

\subsection{All steps combined}

\bn
\item Decompose $n = 3 m + i$, $i \in \{0,1,2\}$
\item Compute the dual representation $ D $ of $ m $.
\item In $ D $, replace every zero with $ Q $ and
  every 1 with $ QM $ and get the word $ W $.
\item Set $ (a, b, c): = (1,0,0) $ and, following the word $ W $
  from left to right, perform the following operations
  modulo $ n $:
\bes
M &:& (a,b,c):= (a+b,b+c,a+b+c) \\
Q &:& (a,b,c):= (a^2 + 2 b c , b^2 + c^2 + 2ab, b^2 + 2ac + 2bc )~.
\ees
\item Finally, calculate
\bes
P_n \mbox{ mod } n = \left\{
\begin{array}{ccc}
3a+2b & \mbox{ for } & i=0 \\
3b+2c & \mbox{ for } & i=1  \\
2a + 2b + 3c   & \mbox{ for } & i=2 ~.
\end{array}
\right.
\ees
\en
For Example we test whether 19 divides $P_{19}$?
\bn
\item $19 = 3 \cdot 6 + 1$, $m=6$, $i=1$
\item Dual representation of 6: $D=110$.
\item $W=QMQMQ$
\item ~~~~~~$\,(a,b,c)=(1,0,0)$\\
  $\stackrel{Q}{\darfo}~(a,b,c)=(1,0,0)$\\
  $\stackrel{M}{\darfo}~(a,b,c)=(1,0,1)$\\
  $\stackrel{Q}{\darfo}~(a,b,c)=(1,1,2)$\\
  $\stackrel{M}{\darfo}~(a,b,c)=(2,3,4)$\\
  $\stackrel{Q}{\darfo}~(a,b,c)=(9,18,11)$
\item $(9,18,11) ~ \stackrel{i=1}{\darfo} ~
3 \cdot 18 + 2 \cdot 11 = 76 \equiv 0 \mbox{ mod } 19$
\en
Thus, we have $19|P_{19}$ and therefore 19 can be a Perrin pseudoprime
or a prime.

\newpage

\subsection{A {\it mathematica}-code for the algorithm}

To deal with large integers we used {\it mathematica}. 
Of course, as an interpretive language it is slower than a compiled code.
But that saved us the development of an own long integer operation package.

The following {\it mathematica}-code was used to check a given
number $n$ on whether $n|P_n$. The code outputs {\tt True} if $n$ is
prime or a Perrin pseudoprime and {\tt False} otherwise. 
We used {\it mathematica11.3} at a
{\tt Intel(R) Core(TM) i5-6500 CPU @ 3.20GHz}.
To check the largest known 1436-digit PPP
(see page \pageref{lang}) takes 0.18 seconds. Checking the largest
Mersenne prime  known in 1982 $2^{86243} - 1$ takes 4 minutes.
Though, at that time the computers were slower.
Today, testing $2^{1398269} - 1$, the 35-th Mersenne prime, found in 1996,
takes a day.
  
\begin{verbatim}
PPP[n_] := (i = Mod[n, 3];
  k = Quotient[n, 3];
  lk = IntegerDigits[k, 2];
  b1 = 1; b2 = 0; b3 = 0;
  Do[ If[ lk[[j]] == 0,
          c1 = b1 * b1 + 2 * b2 * b3;
          c2 = b2 * b2 + b3 * b3 + 2 * b1 * b2;
          c3 = b2 * b2 + 2 * b1 * b3 + 2 * b2 * b3 ,
          a1 = b1 * b1 + 2 * b2 * b3;
          a2 = b2 * b2 + b3 * b3 + 2 * b1 * b2;
          a3 = b2 * b2 + 2 * b1 * b3 + 2 * b2 * b3;
          c1 = a1 + a2;
          c2 = a2 + a3;
          c3 = a1 + a2 + a3];
      b1 = Mod[c1, n]; b2 = Mod[c2, n]; b3 = Mod[c3, n],
   {j, 1, Length[lk]}];
   Which[i == 0, b = 3 * b1 + 2 * b2,
         i == 1, b = 3 * b2 + 2 * b3,
         i == 2, b = 2 * b1 + 2 * b2 + 3 * b3];
   Mod[b, n] == 0)
\end{verbatim}

The Table on \cite{stephan} can be tested with
\begin{verbatim}
ppp = << PPP-new-math;
Do[ If[ Not[ PPP[ ppp[[k1]] ] ] || PrimeQ[ ppp[[k1]] ],
      Print[ ppp[[k1]]," is not a PPP!" ] ], {k1, 1, Length[ppp]}]
\end{verbatim}
Do not forget the semicolon, the list {\tt ppp} is very large.
It runs less than two hours.

\newpage

\section{How to reduce the number of candidates}

It takes many weeks to calculate the 1700 PPP up to $10^{14}$
even with high performance algorithms and computers.
One has to check every number (except a few ones
like mentioned at page \pageref{bauer}
that can be sorted out in advance).
Thus, there is no hope, that one could calculate all
PPPs, say, up to $10^{20}$ in the next years. Moreover,
since they are very rare, if you
take a random $n$, you will ''never'' get a PPP.

So, to calculate more PPPs, one must try to limit the set of
potential candidates.

Dana Jacobsen tested other pseudoprimes, hoping that, for example many
of the Fermat$_2$-PP are also PPPs. And indeed, she found
101994 PPPs up to $18446724258335155361 < 10^{20}$ \cite{dana}.

It turns out that 510 of the 1700 PPPs less than $10^{14}$ are
Fermat$_2$-PP, too.

\subsection{The structure of most of the PPPs}

Let's  have a look at the first PPPs and factorize them:

~

\begin{tabular}{rclcl}
271441 &=& $521 \cdot 521 $ &=& $[ 1(521-1) + 1] \cdot 521 $ \\
904631 &=&  $7 \cdot 13 \cdot 9941 $ \\
16532714 &=&  $2 \cdot 11 \cdot 11 \cdot 53 \cdot 1289 $ && \\
24658561 &=&  $19 \cdot 271 \cdot 4789 $ \\ 
27422714 &=&  $2 \cdot 11 \cdot 11 \cdot 47 \cdot 2411 $ \\
27664033 &=&  $3037 \cdot 9109 $ &=& $[ 3(3037-1) + 1] \cdot 3037$ \\
46672291 &=&  $4831 \cdot 9661 $ &=& $[ 2(4831-1) +1] \cdot 4831 $\\ 
102690901  &=& $5851 \cdot 17551 $ &=& $[ 3(5851-1) +1] \cdot 5851 $\\ 
130944133 &=&  $6607 \cdot 19819 $ &=& $[ 3(6607-1) +1] \cdot 6607 $\\ 
196075949 &=&  $5717 \cdot 34297 $ &=& $[ 6(5717-1) +1] \cdot 5717 $\\ 
214038533 &=&  $8447 \cdot 25339 $ &=& $[ 3(8447-1) +1] \cdot 8447 $\\ 
517697641 &=&  $6311 \cdot 82031 $ &=& $[ 13(6311-1) +1] \cdot 6311 $\\ 
545670533 &=&  $13487 \cdot 40459 $ &=& $[ 3(13487-1) +1] \cdot 13487 $\\
801123451 &=&  $8951 \cdot 89501 $ &=& $[ 10(8951-1) +1] \cdot 8951 $\\ 
855073301 &=&  $16883 \cdot 50647 $ &=& $[ 3(16883-1) +1] \cdot 16883 $\\
903136901 &=&  $17351 \cdot 52051 $ &=& $[ 3(17351-1) +1] \cdot 17351 $\\
970355431 &=&  $22027 \cdot 44053 $ &=& $[ 2(22027-1) +1] \cdot 22027 $\\
\end{tabular}

~

We see that many of them have the structure
$P = [k(p-1) +1] \cdot p$, with some $p \in \P$ and
$k=1,2,3,...$ is a small number.
Clearly, such numbers are never prime. Moreover,
to calculate numbers $P$ in the region of $10^{16}$, it is sufficient to
consider factors $\sim 10^{8}$. Thus, taking into account that
we have 5761455 primes up to $10^{8}$, we get all
pseudoprimes of this structure up to $\sim 10^{16}$ for a given $k$
in half an hour.

This was the starting point of a couple of ideas to
reduce the amount of candidates to be tested.
We list them here in their logical order.

\bn
\item Consider numbers $P=\big[k(p-1)+1\big]p $, ~~$p \in \P$\\
  It was amazing that already $k=3$ and $k=2$ gives more than
  $ 50 \%$ of the 1700 known PPPs up to $10^{14}$.
\item Next, we considered numbers like
  $P=\big[k_1(p-1)+1\big] \big[k_2(p-1)+1\big]$, $p \in \P$; 
  $\mbox{gcd}(k_1,k_2)=1$.   
\item We saw that some PPPs of this structure were overlooked, because
  $p$ must not be prime. Thus, we considered numbers like
  $P=\big[k_1(p-1)+1\big] \big[k_2(p-1)+1\big]$, $p \not\in \P$, $p$ odd.
\item Clearly, the next step were numbers of the form\\
  $P=\big[k_1(p-1)+1\big] \big[k_2(p-1)+1\big] \big[k_3(p-1)+1\big]$
\item and generally $P = \prod_{i=1}^m \big[k_i(p-1)+1\big]$.
  For $m>3$ we get only a few new PPP's.
\en
With this method, we calculated all PPP's with 2 factors
for given $k_i < 100$, with 3 factors for given $k_i < 15$,
and with 4 factors for $k_i < 10$ up to $10^{20}$.
More than $ 95 \%$ of the 1700 known PPPs up to $10^{14}$ have
such a structure. Extrapolating this result, we assume
that we know now  $ 95 \%$ of the PPPs up to $10^{20}$.

It was not possible to find such a PPP with 5 factors for months.

The largest PPPs have about 40 digits.

~

To calculate larger PPPs we used two different methods:
\bi
\item Starting from a PPP with $m$ factors, guess a PPP with $m+1$ factors
  with the same $p$ and some $k_{m+1}$ resulting form the other
  $k_1,...,k_m$. For example, take $k_{m+1}$ as a multiple of the
  least common multiple of the $k_1,...,k_m$. In this way we could find some
  very large PPPs.
\item Do we have to test all odd $p$?  It turns out that
  only a few remainders of $p$ with respect to  23 occur.
  In this way we could find millions of new PPPs up to $10^{24}$.
  \ei
  
\subsection{The remainders of $p$}

Since 23 is the discriminant of the corresponding polynomial of the
Perrin sequence, we look at the 
remainders of $p$ with respect to 23 in more detail.
It turns out that for a given pair $(k_1,k_2)$ we have only a few
remainders instead of 23 possible ones.

For example:
\bi
\item
Take $(k_1,k_2)=(3,1)$, we have the remainders $=(1, 2, 6, 9, 18)$
\item
Take $(k_1,k_2)=(2,1)$, we have the remainders $=(1, 2, 13, 16, 18)$
\ei

The same holds for multiples of 23.
Taking, for example, the number
$23 \cdot 2 \cdot 3 \cdot 5 \cdot 7  \cdot 11 \cdot 13 =690690$.
We have
\bi
\item
For $(k_1,k_2)=(3,1)$ only 14853 remainders (a proportion of 0.0215046),
\item
For $(k_1,k_2)=(2,1)$ only 7425 remainders (a proportion of 0.0107501).
\ei
During our calculation we considered the remainders
with respect to $23 \cdot 2 \cdot 3 = 138$.

Here is a collection of the remainders with respect to 138 for
all pairs $(k_1,k_2)$ with $k_1=5$ and $k_1=7$:

\begin{tabular}{c|c|l}
$k_1$ & $k_2$ & possible remainders with respect to 138 \\ \hline
5 & 1 & 1, 25, 31, 55, 73, 121 \\
5 & 2 & 1, 7, 15, 21, 25, 43, 61, 67, 93, 99, 117, 135 \\
5 & 3 & 1, 9, 25, 43, 55, 63, 75, 93, 109, 117, 121, 135 \\
5 & 4 & 1, 7, 31, 43, 67, 73 \\
7 & 1 & 1, 13, 25, 29, 31, 35, 47, 59, 71, 77, 121, 127 \\
7 & 2 & 1, 13, 25, 67, 97 \\
7 & 3 & 1, 5, 11, 19, 25, 29, 47, 65, 71, 97, 103, 121 \\
7 & 4 & 1, 11, 13, 19, 31, 47, 59, 65, 67, 77, 103, 113 \\
7 & 5 & 1, 25, 31, 67, 121 \\
7 & 6 & 1, 5, 13, 29, 47, 59, 67, 79, 97, 113, 121, 125 \\
\end{tabular}

~

%\rb{Wie hingen die zusammen? Genauere cojectures}

These remainders were found experimentally. For a given
pair $(k_1,k_2)$ we calculated some PPPs for any odd $p$,
enough to be sure about the possible remainders. Having obtained
these, we test the following $p$ only with these remainders.
That resulted in a strong speed-up.

Unfortunately, we have no idea how the remainders can
be calculated in advance. We think this is an interesting
problem for specialists, for example, in Carmichael numbers.

~

For PPPs with 3 factors we observed the following interesting
experimental result:

Fix a pair $(k_1,k_2)$ with $gcd(k_1,k_2)=1$ and let be
$R(k_1,k_2)$ the set of remainders of $p$. Then, the set of
remainders $R(k_1,k_2,k_3)$ for a PPP with 3 factors is
\bes
R(k_1,k_2,k_3) = R(k_1,k_2) \cap  R(k_1,k_3) \cap  R(k_2,k_3)
\ees
Thus, the number of possible remainders decreases with the number of factors.

A similar result holds for  PPPs with more than 3 factors.
Again, we do not know how to prove this.

~

The remainder 1  with respect to multilpes of 23 contains in any set
of remainders for any $(k_i)$.

%\newpage

\section{Numerical results}

\subsection{The state of the art}

A current overview can be found in N.J.A. Sloanes famous OEIS 
(On-Line Encyclopedia of Integer Sequences) \cite{oeis}.

By now, all PPPs -- 1700 -- up to $10^{14}$ are known.
Since we have 3204941750802 primes up to $10^{14}$,
using the Perrin prime test, a PPP occurs with probability
$W(10^{14})=5.30431 10^{-10}$.
Thus, to check whether a given number less than $10^{14}$ is prime you can use
the Perrin test and -- if it is true -- look at the table whether it is one
of the 1700 PPPs. If not, it is prime.

The following table shows the probability 
$W(n)$ up to $n=10^{14}$. We used \cite{primes} for the numbers of primes.

~

\begin{tabular}{l|r|r|r}
  $n$ & PPPs & primes & probability $W(n)$ \\ \hline
 $ 10^8$   &   7 &        5761455 & $  1.21497 *10^{-6}~$ \\
$10^9 $    &  17 &       50847534 & $  3.34333 *10^{-7}~$ \\
$10^{10}$  &  42 &      455052511 & $  9.22970  *10^{-8}~$ \\
$10^{11}$  & 116 &     4118054813 & $  2.84115 *10^{-8}~$ \\
$10^{12}$  & 285 &    37607912018 & $ 7.57819 *10^{-9}~$ \\
$10^{13}$  & 649 &   346065536839 & $ 1.87537 *10^{-9}~$ \\
$10^{14}$  & 1700 & 3204941750802 & $  5.30431 *10^{-10}$  \\
  \end{tabular}

  ~

  \subsection{Our results}

  We calculated {\bf 9261931} (by December 2019)
  PPPs that an be found in the database \cite{stephan}.
  (Note, that the database is updated from time to time.)
  
  We tried to find all PPPs up to $10^{20}$ 
  and all with 2 factors and $(k_1,k_2)=(3,1)$ and $(k_1,k_2)=(2,1)$
  up to $10^{22}$. Of course there is a by-catch of many PPPs up to $10^{30}$.
  
  Moreover, we tried to find some very large ones using two methods:

  At first, we constructed PPPs with $m+1$ factors starting from a known ones
  with $m$ factors.

  Second, knowing that 1 is always a remainder with respect to multilpes of 23
  for all $p$, we tested numbers of the form $n = p \cdot \big( k(p-1)+1 \big)$.
  with $k=2,3$ and
  $p = 23 \cdot 2 \cdot 3 \cdot 5 \cdot 7 \cdot 11 \cdots$ a multiple
  of 23 and the first primes
  This yields very large PPPs, for example the one on page \pageref{lang}.

 \subsubsection{Almost all PPPs}

\begin{minipage}[c]{8cm}
Having a look at the table above, we see that
$\log W(n)$ behaves largely linearly. We extrapolate this and 
expect the following numbers of PPPs. The numbers up to $10^{20}$ are
``almost all'', the  numbers up to $10^{22}$ are ``more than a half''
of all PPPs.  
\unitlength=1cm
\begin{picture}(8.0,6.0)\thicklines
  \put(0.0,0.0){\includegraphics[width=8.0cm]{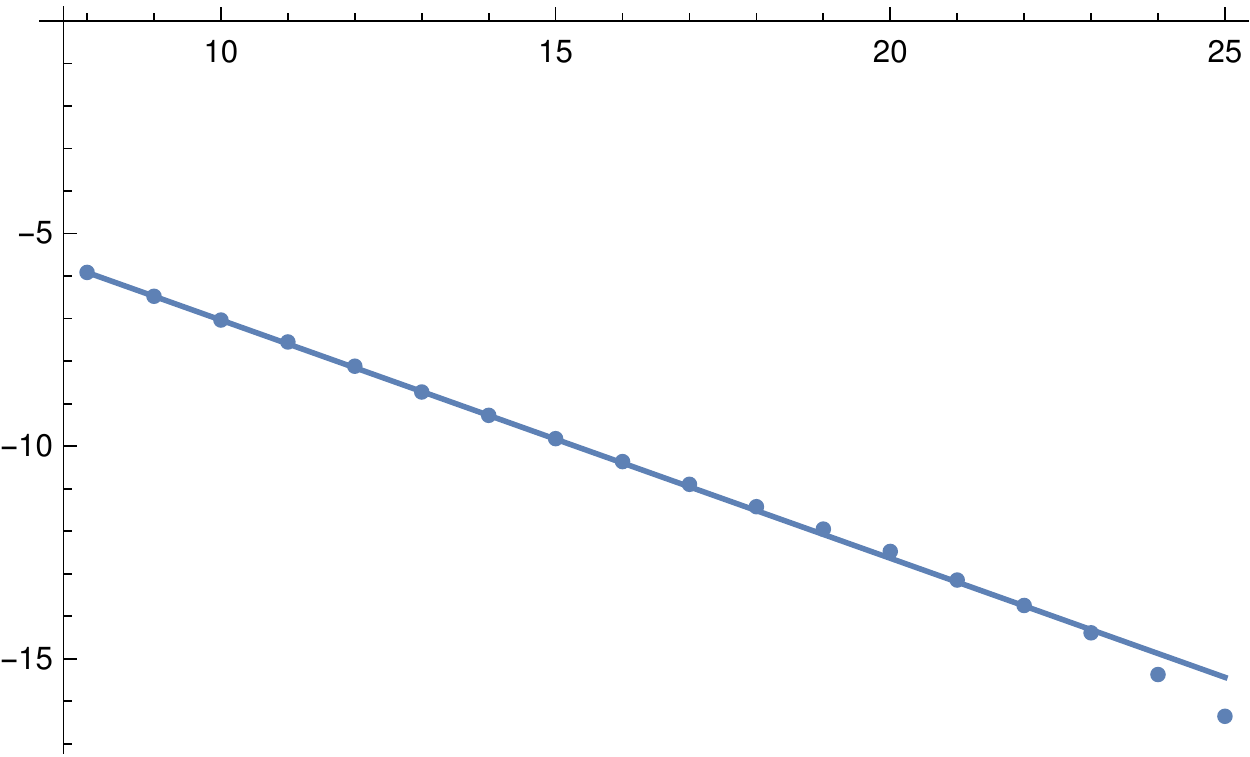}}
 \put(7,5.3){\makebox(0.0,0.0){$\log_{10} n$}}
 \put(2.0,0.5){\makebox(0.0,0.0){$\log_{10} W(n)$}}
\end{picture}

\end{minipage}\hfill
\begin{minipage}[c]{7cm}
\begin{tabular}{l|r|r}
  $n$ & expected PPPs  & founded PPPs \\ \hline
 $10^{15}$ &        4360  & 4409 \\
 $10^{16}$ &       11236  & 11972 \\
 $10^{17}$ &       29076  & 33045 \\
 $10^{18}$ &       75520  & 93001 \\
 $10^{19}$ &      196790 & 262236 \\
 $10^{20}$ &      514287 & 742759 \\
 $10^{21}$ &     1347560  & 1502883 \\
 $10^{22}$ &     3539332  & 3615622 \\
 $10^{23}$ &     9316050  & 7870747 \\
 $10^{24}$ &    24569601  & 7874995 \\
 $10^{25}$ &    64915566  & 7879187 \\
 $10^{26}$ &   171799266  & 7885930 \\
 $10^{27}$ &   455365341  & 7898184 \\
 $10^{28}$ &  1208691635  & 7920907 \\
 $10^{29}$ &  3212505576  & 7964655 \\
 $10^{30}$ &  8548808804  & 8049285 \\
\end{tabular}
\end{minipage}

  \newpage

 \subsubsection{Huge PPPs}

{\bf Collected by factors:} We found
\bi
\item 1 PPP with 14 factors.
\item 13 PPPs with 13 factors.
\item 64 PPPs with 12 factors.
\item 113 PPPs with 11 factors.
\item 176 PPPs with 10 factors.
\item 481 PPPs with 9 factors.
\item 1054 PPPs with 8 factors.
\item 2591 PPPs with 7 factors.
\item 7159 PPPs with 6 factors.
\item 29529 PPPs with 5 factors.
\ei
%This was only possible by restricting to a reduced set of remainders.

~

{\bf Collected by digits:}  We found
\bi
\item $\sim 4000$ PPPs with more than 80 decimal digits
\item $\sim 1600$ PPPs with more than 100 decimal digits
\item $36$ PPPs with more than 500 decimal digits
\item $6$ PPPs with more than 1000 decimal digits
\item The largest PPP has 3101 digits. Here it is:
  {\tiny
\bes
\hspace*{-5em}
&&2182001064371918934845924375655593970781204553917566660863280384747887616030277480053172205785183353188400\\[-0.3em]\hspace*{-5em}
&&4126146210865090197070653868880189559625867459754727073713090924616711853613422828119114381617102058517546\\[-0.3em]\hspace*{-5em}
&&8653751496284195684000100419880283999039015488001095163810247785156033211221423472140681188918922518742770\\[-0.3em]\hspace*{-5em}
&&0398996872031544022682029689624783660853880129295123479427747681652039459239579760489615206781614707161883\\[-0.3em]\hspace*{-5em}
&&9138537548347177754556329233097993446947475927879595917904730731452471057039913228447069819231974147528469\\[-0.3em]\hspace*{-5em}
&&7693616171472498459173243671532936165356214403017220481995761095314765972379574827945192124085559691984391\\[-0.3em]\hspace*{-5em}
&&8008661242667729379149221402733564699474653803584334247108722459604844155931040562979301921938928545995807\\[-0.3em]\hspace*{-5em}
&&4207926519074011909871332364749649617141024864366985374867133374038568149858039921667907016960062202008122\\[-0.3em]\hspace*{-5em}
&&9182067899216118132468035588845067378082718617393902077009092862097562284582389695785019716348129717066692\\[-0.3em]\hspace*{-5em}
&&0783325505675383114442119375756418942531432620905077133117297177064802424569877645651274316923030865339422\\[-0.3em]\hspace*{-5em}
&&6661109617675061215430499075868542147459797368102792867066735398199032669816585264700339738266181367925685\\[-0.3em]\hspace*{-5em}
&&9183901438799475057989326512787989244219170992158347364160368593405317157057039942593979747214483064168779\\[-0.3em]\hspace*{-5em}
&&3723363454025576455261406877507795872082604992320378872519383088242811076665512015332176716276340248257164\\[-0.3em]\hspace*{-5em}
&&6729443535184738262902790223792682930259972646770066028255813046639125749771256788743514165965139691554159\\[-0.3em]\hspace*{-5em}
&&3353592560965482315120431456622925845399082336306306234166863238919515156950417488352070194395498058003429\\[-0.3em]\hspace*{-5em}
&&2609689928226091668646468088635185719074533550653987615133601688385577315810376211381436151897390975873498\\[-0.3em]\hspace*{-5em}
&&9194775781036920280653165835092015711042583063595692979056408307560965084104645943087850367750725513620664\\[-0.3em]\hspace*{-5em}
&&9589379996405514942415050679736879467176251813294056719410189773891939434281262409431885675830573414891359\\[-0.3em]\hspace*{-5em}
&&7068260880092249389030829673092944201188379579217564895495418187279934349004962876837044167260718567772046\\[-0.3em]\hspace*{-5em}
&&7521150708667751876125544569499435754902575963129390715770989789849330459963345038762428879760367628428833\\[-0.3em]\hspace*{-5em}
&&7083464467875818139474195085529183097604033933360012552535245232509900842279440109453302234497800743667133\\[-0.3em]\hspace*{-5em}
&&2290093368659872164696682455863309852162786109791145473780233128398296687924256984146263917624053810047106\\[-0.3em]\hspace*{-5em}
&&9132240022024999815261877155099328326233538506570393468310793807821070234336347574184496483617336881484518\\[-0.3em]\hspace*{-5em}
&&9783926914876429525603769119738558257277589955344693025872664161546365759997766592490233729898293133230624\\[-0.3em]\hspace*{-5em}
&&4301770299046097662381531807593304842496115443710755824125123112656492287865978030693101114925766670096297\\[-0.3em]\hspace*{-5em}
&&4043457120990040352730767662860730019992114778921176312285224644592166173374663104973515972020108030670776\\[-0.3em]\hspace*{-5em}
&&0538966132268173354370805800388713443173563909282726774947019900416544732774260586167631835100825092596248\\[-0.3em]\hspace*{-5em}
&&4432038054992189389231847184387110810917603905274409490013362690801082371949435532760468825732391337145460\\[-0.3em]\hspace*{-5em}
&&6507376646884319008228201004154992411941387896249068825523566890040592991334780411481021215235342677940980\\[-0.3em]\hspace*{-5em}
&&162869702039217052132582551\\
\ees
}
\ei
\label{lang}

%\rb{Factorising is left to the reader}

\subsubsection{Some more information}

\bi
\item Our method found 1647 out of the known 1700 up to $10^{14}$. Thus,
  53 or $\sim 3\%$ left. We call them ``sporadic PPPs''.
\item Dana Jacobsen's list of 101994 PPPs contains 699 that we could not
  find with our method.
\item We found 742759 PPPs up to $10^{20}$.
  If these compile $97\%$ of all PPPs, then 22972 sporadic ones are left.
\item Among the the first 10000 Carmichael numbers (taken from \cite{car})
  there are 16 PPPs:
  \bes
C_{1353} &=& 7045248121 = 821 * 1231 * 6971 \eyy
(2 (411-1) + 1) * (3 (411-1) + 1) * (17 (411-1) + 1)\\
C_{1375} &=& 7279379941 = 211 * 3571 * 9661\\
C_{2142} &=& 24306384961 = 19 * 53 * 79 * 89 * 3433\\
C_{2652} &=& 43234580143 = 223 * 5107 * 37963\\
C_{2837} &=& 52437986833 = 23 * 463 * 1453 * 3389\\
C_{2988} &=& 60518537641 = 23 * 89 * 991 * 29833\\
C_{3336} &=& 80829302401 = 89 * 199 * 463 * 9857\\
C_{3855} &=& 118805562613 = 829 * 9109 * 15733\\
C_{4125} &=& 144377609419 = 1319 * 9227 * 11863\\
C_{4322} &=& 165321688501 = 101 * 271 * 691 * 8741\\
C_{4342} &=& 167385219121 = 83 * 6971 * 289297\\
C_{5046} &=& 254302215553 = 307 * 3673 * 225523\\
C_{5731} &=& 364573433665 = 5 * 7 * 23 * 37 * 997 * 12277\\
C_{6743} &=& 575687567521 = 11 * 19 * 79 * 137 * 307 * 829\\
C_{6810} &=& 588909469501 = 1871 * 16831 * 18701 \eyy
1871*( 9(1871-1)+1)*( 10(1871-1)+1)\\
C_{7057} &=& 652270080001 = 3361 * 9241 * 21001
\ees
Some of them, namely,
$C_{2142}$, $C_{2837}$, $C_{3336}$, $C_{4342}$, $C_{5731}$, 
$C_{6743}$ and $C_{7057}$ we could not find with our method.

Note, that
$C_{7057} = 
(  4*(841 - 1) + 1) * (  11*(841 - 1) + 1) * (  25*(841 - 1) + 1) $
with $841 = 19^2$. We could not find it, since we
restrict ourself to $k_i \leq 15$ for numbers with 3 factors.
\ei

\subsubsection{Some conjectures}

During the calculations, we were led to the following 
conjectures. We invite everyone to think about the proofs.

\bi
\item Almost all PPPs have the structure
  $P = \prod_{i=1}^m \big[k_i(p-1)+1\big]$
\item There are infinitely many of such type.
\item The $p$ has few remainders with respect to multiples of 23.
  They can be calculated theoretically in advance.
\item If $\prod_{i=1}^m \big[k_i(p-1)+1\big]$ is a PPP, then with
``high'' probability \\[0.4em]
$\prod_{i=1}^{m+1} \big[k_i(p-1)+1\big]$ is a PPP
with $k_{m+1} = c k_{m}$. In such a way you can 
construct large PPPs.
\item The set of remainders (with respect to multiples of 23) of $p$
  corresponding to given $k_i$ with 3 (or more) factors are the
  intersection of the sets of remainders corresponding to fewer $k_i$,
  requiring $gcd(k_i,k_j)=1$.
\item There are a particularly large number of PPPs if the $k_i$
  are prime, pairwise.
\item If for some $p$ the number with $\{k_2 \cdot k_3,k_2,k_3\}$ is
  a PPP then so is the number with $\{k_2,k_3\}$. 
\ei

\section{Other promising polynomials for pseudoprimes}

We tested polynomials of degree 3 and 4
with integer coefficients $a_i$ with $|a_i| \leq 20$.
Every corresponding sequences we tested for pseudoprimes up to $10^9$.
For polynomials of third order the Perrin sequence is indeed the rarest.

For polynomials of fourth order we find two polynomials without any
pseudoprimes up to $10^9$ at all. Here they are:
\bes
Q(x) &=& - x^4 + x^3 - 17 x^2 + 0 x + 5 \\
R(x) &=& - x^4 + 11 x^3 + x^2 - 12 x + 14 
\ees
We have for $Q(x)$ the corresponding sequence
\bes
q_n &=& q_{n-1} - 17 q_{n-2} + 5 q_{n-4} \\
q_0 &=& 4 \\
q_1 &=& 1 \\
q_2 &=& -33 \\
q_3 &=& -50 
\ees
and the testing rule $n \in \P$ $\darfo$ $n|(q_n-1)$.

For $R(x)$ the sequence is
\bes
r_n &=& 11 r_{n-1} + r_{n-2} - 12 r_{n-3} + 14 r_{n-4} \\
r_0 &=& 4 \\
r_1 &=& 11 \\
r_2 &=& 123 \\
r_3 &=& 1328 
\ees
and the testing rule is $n \in \P$ $\darfo$ $n|(r_n-11^n)$.

To avoid the term $11^n$, it is better to consider
\bes
G(x) = Q(x) (x-11) = - x^5 + 22 x^4 - 120 x^3 - 23 x^2  + 146 x -154
\ees
instead of $R(x)$. This corresponds to the 5-th oder sequence
\bes
g_n &=& 22 g_{n-1} - 120 g_{n-2} -23 g_{n-3} + 146 g_{n-4} - 154 g_{n-5} \\
g_0 &=& 3 \\
g_1 &=& 0 \\
g_2 &=& 2 \\
g_3 &=& -3 \\
g_4 &=& 14 
\ees
with the testing rule $n \in \P$ $\darfo$ $n|g_n$.

%\newpage


\begin{thebibliography}{00}
\bibitem{granv1} W. R. Alford, A. Granville, C. Pomerance,
There are Infinitely Many Carmichael Numbers, Ann. Math. 139, 703-722, 1994.
\bibitem{alg} F. Bauern\"oppel, private communication
\bibitem{gauss}C. F. Gauss, Article 329 of
  Disquisitiones Arithmeticae (1801)
\bibitem{grant} J. Grantham: There are infinitely many Perrin pseudoprimes.
Journal of Number Theory. 130, Nr. 5, 2010, S. 1117-1128  
\bibitem{dana} D. Jacobsen, \verb+http://ntheory.org/pseudoprimes.html+
\bibitem{stephan}
  H. Stephan, Perrin pseudoprimes. Data Sets,
  Weierstrass Institute Berlin (2019),
  \verb+http://doi.org/10.20347/WIAS.DATA.4+
\bibitem{wiki1} \verb+https://en.wikipedia.org/wiki/Perrin_number+
\bibitem{oeis} \verb#https://oeis.org/search?q=perrin+pseudoprimes#
\bibitem{car} \verb+https://oeis.org/A002997/b002997.txt+
\bibitem{primes} \verb+https://primes.utm.edu/howmany.html+
\end{thebibliography}
\end{document}